\documentclass[10pt]{article}

\usepackage{graphicx}%
\usepackage{amsmath,amssymb,amsthm,upref,bm}%

\newtheorem{theorem}{Theorem}%
\newtheorem{remark}{Remark}

\begin{document}

\baselineskip=4.4mm

\makeatletter

\newcommand{\E}{\mathrm{e}\kern0.2pt} 
\newcommand{\D}{\mathrm{d}\kern0.2pt}
\newcommand{\RR}{\mathbb{R}}
\newcommand{\CC}{\mathbb{C}}%
\newcommand{\ii}{\kern0.05em\mathrm{i}\kern0.05em}

\renewcommand{\Re}{\mathrm{Re}} 
\renewcommand{\Im}{\mathrm{Im}}

\def\bottomfraction{0.9}

\title{\bf ``Potato kugel'' for nuclear forces and \\ a small one for acoustic
waves}

\author{Nikolay Kuznetsov}

\date{ }

\maketitle

\vspace{-8mm}

\begin{center}
Laboratory for Mathematical Modelling of Wave Phenomena, \\ Institute for Problems
in Mechanical Engineering, Russian Academy of Sciences, \\ V.O., Bol'shoy pr. 61, St
Petersburg 199178, Russian Federation \\ E-mail: nikolay.g.kuznetsov@gmail.com
\end{center}

\begin{abstract}
\noindent The ``potato kugel'' theorem of Aharonov, Schiffer and Zalcman, which
concerns an inverse property of harmonic functions, is extended to the settings of
the modified Helmholtz and Helmholtz equations that describe nuclear forces and
acoustic waves respectively.
\end{abstract}

\setcounter{equation}{0}

\vspace{-4mm}

\section{Introduction}

Analytic characterization of balls in the Euclidean space $\RR^m$ by means of
harmonic functions has a long history; it started in the 1960s, in the pioneering
notes \cite{E}, \cite{ES}, and shortly after\-wards Kuran \cite{K} obtained the
following general result:

\vspace{2mm}

\noindent {\bf Theorem K.} {\it Let $D$ be a domain (= connected open set) of
finite (Lebesgue) measure in the Euclidean space $\RR^m$ where $m \geq 2$. Suppose
that there exists a point $P_0$ in $D$ such that, for every function $h$ harmonic in
$D$ and integrable over $D$, the volume mean of $h$ over $D$ equals $h (P_0)$. Then
$D$ is an open ball (disk when $m=2$) centred at $P_0$.}

\vspace{2mm}

Presumably, the paper \cite{HN} was the first one in which this theorem was referred
to as the property of harmonic functions inverse to the mean value identity for
balls. The term became widely accepted. A slight modification of Kuran's
considerations shows that his theorem is valid even if~$D$ is disconnected; see the
survey article \cite{NV}, p.~377, which also contains some improvements of Kuran's
theorem, and a discussion of its applications and of possible similar results
involving certain averages over $\partial D$, when $D$ is a bounded domain.

Another approach to harmonic characterization of balls was developed by Aharonov,
Schiffer and Zalcman \cite{ASZ} (the origin of a rather unusual title of their paper
is explained in Zalcman's comment; see \cite{S}, p.~497). They proved the following:

\vspace{2mm}

\noindent {\bf Theorem ASZ.} {\it Let $D \subset \RR^m$, $m \geq 3$, be a bounded
open set. If the equality
\begin{equation*}
\int_D |y - x|^{2-m} \D y = \frac{a}{|x|^{m-2}} + b \label{asz}
\end{equation*}
holds with suitable real constants $a$ and $b$ for every $x \in \RR^m \setminus D$,
then $D$ is an open ball centred at the origin, $a = |D|$ and $b=0$.}

\vspace{2mm}

Here and below $|D|$ stands for the volume of $D$ (area if $D \subset \RR^2$). Since
$|y - x|^{-1}$ is a fundamental solution of the Laplace equation for $m=3$, this
theorem answers in the affirmative the following question posed to the authors of
\cite{ASZ}:
\begin{quote}
Let $D$ be a homogeneous, compact, connected ``potato'' in space, which
gravitationally attracts each point outside it as if all its mass were concentrated
at a single point. Does this guarantee that $D$ is a ball centred at this point?
\end{quote}
Later, Lanconelli \cite{L} extended Theorem ASZ to the sub-Laplacian setting, and in
the recent article \cite{CL}, an improvement of Theorem ASZ was obtained. It relaxes
the original restriction imposed on $D$, and a more natural identity is used to
guarantee that $D$ is a ball:

\vspace{2mm}

\noindent {\bf Theorem ASZ$'$} [Cupini, Lanconelli]. {\it Let $D \subset \RR^m$, $m
\geq 3$, be an open set such that $|D| < \infty$. If for some $x_0 \in D$ the
identity
\begin{equation}
|D|^{-1} \int_D |y - x|^{2-m} \D y = |x_0 - x|^{2-m} \label{asz'}
\end{equation}
holds for every $x \in \RR^m \setminus D$, then $D$ is an open ball centred at
$x_0$.}

\vspace{2mm}

The recent survey \cite{Ku4} complements \cite{NV} providing coverage for other
results in this area; in particular, various characterizations of balls via
harmonic functions as well as characterizations of other domains (strips, annuli
etc.) are given. Moreover, a characterization of balls via solutions to the modified
Helm\-holtz equation
\begin{equation}
\nabla^2 u - \mu^2 u = 0 , \quad \mu \in \RR \setminus \{0\} \, . \label{MHh}
\end{equation}
is considered. Here and below, $\nabla = (\partial_1, \dots , \partial_m)$,
$\partial_i = \partial / \partial x_i$, denotes the gradient operator. The most
important application of this equation is in the theory of nuclear forces; it was
developed by Yukawa in his Nobel Prize winning paper \cite{Y}. See also \cite{D},
where \eqref{MHh} is referred to as the Yukawa equation and its solutions are called
{\it panharmonic functions}---abbreviation used below. Since the mentioned
characterization of balls obtained in \cite{Ku} is closely related to one of the
results proved in this paper, we begin with its formulation, but introduce some
notation before that.

Let $x = (x_1, \dots, x_m)$ be a point in $\RR^m$, $m \geq 2$, by $B_r (x) = \{ y
\in \RR^m : |y-x| < r \}$ we denote the open ball of radius $r$ centred at $x$ (just
$B_r$, if centred at the origin). The ball is called admissible with respect to a
domain $D \subset \RR^m$ provided $\overline{B_r (x)} \subset D$. If $D$ has finite
Lebesgue measure and a function $f$ is integrable over $D$, then
\[ M (f, D) = \frac{1}{|D|} \int_{D} f (x) \, \D x
\]
is its volume mean value over $D$. The volume of $B_r$ is $|B_r| = \omega_m r^m$,
where
\[ \omega_m = 2 \, \pi^{m/2} / [m \, \Gamma (m/2)]
\]
is the volume of the unit ball; here $\Gamma$ denotes the Gamma function. A dilated
copy of a domain $D$ is $D_r = D \cup \left[ \cup_{x \in \partial D} B_r (x)
\right]$. Thus, the distance from $\partial D_r$ to $D$ is equal to $r$.

Now, we formulate an analogue of Theorem~K, which was proved in~\cite{Ku}.

\begin{theorem}
Let $D \subset \RR^m$, $m \geq 2$, be a bounded domain, whose complement is
connected, and let $r > 0$ be such that $|B_r| = |D|$. If for a point $x_0 \in D$
and some $\mu > 0$ the identity
\begin{equation*}
u (x_0) \, a_m^+ (\mu r) = M (u, D) \, , \ \ \mbox{where} \ \ a_m^+ (t) = \Gamma
\left( \frac{m}{2} + 1 \right) \frac{I_{m/2} (t)}{(t / 2)^{m/2}} \, , \label{CR}
\end{equation*}
holds for every positive function $u$ panharmonic in $D_r$, then $D = B_r (x_0)$.
\end{theorem}

As usual, $I_\nu$ denotes the modified Bessel function of order $\nu$. Its
well-known properties (see \cite{Wa}, pp.~79, 80) imply that $a_m^+ (t)$ increases
monotonically for $t \in [0, \infty)$ from $a_m^+ (0) = 1$ to infinity. In the
three-dimensional case related to nuclear forces, $a_3^+ (t) = \sqrt{2 \pi} \, I_{3/2}
(t) / t^{3/2}$.

Clearly, Theorem 1 is inverse of {\it the $m$-dimensional mean value property
\begin{equation}
a_m^+ (\mu r) \, u (x) = M (u, B_r (x)) \, , \label{MM+}
\end{equation}
which holds for every admissible ball $B_r (x)$ provided $u \in C^2 (D) \cap L^1
(D)$ solves \eqref{MHh} in $D$.} Identity \eqref{MM+} was recently obtained by the
author \cite{Ku2}, but it was a surprise to discover that only mean value formulae
for spheres and circles were known earlier for panharmonic functions. As early as
1896, C.~Neumann \cite{NC}, Ch.~9, Sect.~3, derived the following formula for
spheres in $\RR^3$
\begin{equation*}
a_1^+ (\mu r) \, u (x) = M (u, \partial B_r (x)) \, , \label{CN}
\end{equation*}
where $a_1^+ (\mu r) = \sinh (\mu r) / (\mu r)$. Much later, Duffin \cite{D},
pp.~111-112, independently rediscovered the same proof, but in $\RR^2$ with $a_1^+
(t)$ replaced by $a_0^+ (t) = I_0 (t)$.

In order to reformulate the question quoted above for nuclear setting we recall
that Yukawa \cite{Y}, p.~49, described a source of nuclear force located at $y \in
\RR^3$ with the help of the following potential
\begin{equation}
E_\mu^- (x, y) = \frac{\exp \{- \mu |x-y|\}}{|x-y|} \, , \quad \mu > 0 , \quad x \in
\RR^3 \setminus \{ y \} , \label{E-}
\end{equation}
which is a nonnegative fundamental solution of equation \eqref{MHh} decaying rapidly
with the distance. Another fundamental solution of this equation grows with the
distance, namely:
\begin{equation}
E_\mu^+ (x, y) = \frac{\exp \{\mu |x-y|\}}{|x-y|} \, , \quad \mu > 0 , \quad x \in
\RR^3 \setminus \{ y \} . \label{E+}
\end{equation}
The existence of two linearly independent fundamental solutions distinguishes
\eqref{MHh} from the Laplace equation.

For every $r>0$ and arbitrary $x_0 \in \RR^3$, these fundamental solutions define
two families of integrable panharmonic functions
\[ B_r (x_0) \ni y \mapsto  E_\mu^\pm (y, x) \ \ \mbox{parametrised by} \ x \in \RR^3 
\setminus B_r (x_0) \, .
\]
The mean value property \eqref{MM+} yields that
\begin{equation}
a_3^+ (\mu r) \, E_\mu^\pm (x, x_0) = M (E_\mu^\pm (\cdot, x), B_r (x_0)) \label{MME+}
\end{equation}
for every element of these families. This identity is analogous to \eqref{asz'} with
$D$ changed to $B_r (x_0)$, but involves the factor $a_3^+ (\mu r) > 1$ on the
left-hand side. Now, in view of Theorem~ASZ$'$ and identity \eqref{MME+}, it is
natural to expect that the following assertion is true.

\begin{theorem}
Let ``potato'' occupy a bounded domain $D \subset \RR^3$, whose complement is
connected, and let $r$ be such that $|B_r| = |D|$. If for some $x_0 \in D$ the mean
value identity
\begin{equation}
a_3^+ (\mu r) \, E_\mu^\pm (x, x_0) = M (E_\mu^\pm (\cdot, x), D) \label{MMED}
\end{equation}
holds for every $x \notin D$ and for each fundamental solution of equation
\eqref{MHh}, then $D = B_r (x_0)$.
\end{theorem}

It occurs that a characterisation of balls analogous to Theorem~1 is valid for
meta\-harmonic functions; the term is just an abbreviation for `solution to the
Helmholtz equation'
\begin{equation}
\nabla^2 u + \lambda^2 u = 0 , \quad \lambda \in \RR \setminus \{0\} \, .
\label{Hh}
\end{equation}
I.~N. Vekua introduced it in 1943, in his still widely cited article, which was also
published as Appendix~2 to the monograph \cite{Ve2}. In the following assertion
proved recently, $J_\nu$ is the Bessel function of order $\nu$, whereas its $n$th
positive zero is denoted by $j_{\nu,n}$.

\begin{theorem} [\cite{Ku1}]
Let $D \subset \RR^m$, $m \geq 2$, be a bounded domain, whose complement is
connected, and let $r > 0$ be such that $|B_r| = |D|$. Suppose that there exists a
point $x_0 \in D$ such that for some $\lambda > 0$ the identity
\begin{equation*}
u (x_0) \, a_m^- (\lambda r) = M (u, D) \, , \ \ \mbox{where} \ \ a_m^- (t) = \Gamma
\left( \frac{m}{2} + 1 \right) \frac{J_{m/2} (t)}{(t / 2)^{m/2}} \, , \label{new}
\end{equation*}
holds for every function $u$ metaharmonic in $D_r$. If also
\begin{equation}
D \subset B_{r_0} (x_0) , \ \ where \ \ \lambda r_0 = j_{m/2,1} \, ,
\label{exam}
\end{equation}
then $D = B_r (x_0)$.
\end{theorem}

\begin{remark}
{\rm For a fixed $\lambda > 0$, the assertion is applicable only to domains, whose
volume is less than or equal to $|B_{r_0}|$, where $\lambda r_0 = j_{m/2,1}$.
Indeed, every such domain must lie within a ball of radius $r_0$, and this
distinguishes the last theorem from Theorem~1, imposing no restriction on the
domain's volume.

The reason for restriction \eqref{exam} is as follows. The function $a_{m-2}^- (t)$
used in the proof of Theorem~3 oscillates about zero, and so only an interval near
the origin, where $a_{m-2}^-$ decreases monotonically, can be used. On the opposite,
$a_{m-2}^+ (t)$ used in the similar proof of Theorem~1 is greater than or equal to
unity and increases monotonically.}
\end{remark}

Like Theorem 1, the last theorem is inverse of {\it the $m$-dimensional mean value
property
\begin{equation}
a_m^- (\lambda r) \, u (x) = M (u, B_r (x)) \, , \label{MM-}
\end{equation}
which holds for every admissible ball $B_r (x)$ provided $u \in C^2 (D) \cap L^1
(D)$ solves \eqref{Hh} in $D$.} Identity \eqref{MM-} was also obtained by the author
\cite{Ku2} only recently.

The existence of two linearly independent fundamental solutions is another feature
common to equations \eqref{MHh} and \eqref{Hh}, but unlike \eqref{E-} and \eqref{E+}
the solutions of \eqref{Hh}, namely,
\begin{equation}
E_\lambda^\pm (x, y) = \frac{\exp \{\pm \ii \lambda |x-y|\}}{|x-y|} \, , \quad
\lambda > 0 , \quad x \in \RR^3 \setminus \{ y \} , \label{Epm}
\end{equation}
are complex-valued, thus allowing to describe outgoing and incoming acoustic waves
in the time domain; see \cite{Ve2}, Appendix~2.

As in the panharmonic case, for every $r>0$ and arbitrary $x_0 \in \RR^3$, the
fundamental solutions \eqref{Epm} define two families of integrable metaharmonic
functions
\[ B_r (x_0) \ni y \mapsto  E_\lambda^\pm (y, x) \ \ \mbox{parametrised by} \ x \in
\RR^3 \setminus B_r (x_0) .
\]
The mean value property \eqref{MM-} yields that
\begin{equation}
a_3^- (\lambda r) \, E_\lambda^\pm (x, x_0) = M (E_\lambda^\pm (\cdot, x), B_r
(x_0)) \label{MME-} 
\end{equation}
for every element of these families.  This identity is analogous to \eqref{asz'}
with $D$ changed to $B_r (x_0)$, but involves the factor $a_3^- (\lambda r) < 1$ on
the left-hand side, which is positive on the interval $(0, j_{3/2,1})$ and then
changes sign. In view of Theorem~ASZ$'$ and identity \eqref{MME-}, it is natural to
expect that the following assertion is true.

\begin{theorem}
Let ``potato'' occupy a bounded domain $D \subset \RR^3$, whose complement is
connected, and let $r$ be such that $|B_r| = |D|$. If for some $x_0 \in D$ condition
\eqref{exam} is fulfilled for $m=3$ and some $\lambda > 0$, and the mean value
identity
\begin{equation}
a_3^- (\lambda r) \, E_\lambda^\pm (x, x_0) = M (E_\lambda^\pm (\cdot, x), D) 
\label{MMEDH}
\end{equation}
holds for every $x \notin D$ and for each fundamental solution of equation
\eqref{Hh}, then $D = B_r (x_0)$.
\end{theorem} 

\begin{remark}
{\rm In the acoustical case with a fixed wave number $\lambda > 0$, a restriction
on the size of ``potato'' $D$ is imposed by condition \eqref{exam}. Namely, the
volume of $D$ must be less than or equal to $|B_{r_0}|$, where $\lambda r_0 =
j_{3/2,1}$; moreover, every such domain must lie within a ball of radius $r_0$.}
\end{remark}

\section{Proof of Theorems 2 and 4}

\begin{proof}[Proof of Theorem 2]
Since $E_\mu^+ (x, x_0)$ and $E_\mu^- (x, x_0)$ satisfy \eqref{MMED} for every $x
\notin D$, the same is true for every linear combination of these fundamental
solutions. In particular,
\[ |D| \, a_3^+ (\mu r) \frac{\sinh (\mu |x - x_0|)}{\mu |x - x_0|} =
\int_D \frac{\sinh (\mu |x - y|)}{\mu |x - y|} \, \D y \ \ \mbox{for every} \ x
\notin D .
\]
Moreover, the identity is valid throughout $\RR^3$, because we have real-analytic
functions of $x$ on both sides (a consequence of the analyticity of $z^{-1} \sinh
z$), and so we substitute $x=x_0$, thus obtaining
\[ |D| \, a_3^+ (\mu r) = \int_D \frac{\sinh (\mu |x_0 - y|)}
{\mu |x_0 - y|} \, \D y \, .
\]
Let us relocate, without loss of generality, the domain $D$ so that $x_0$ coincides
with the origin, which simplifies the identity to
\begin{equation}
|D| \, a_3^+ (\mu r) = \int_D U_+ (y) \, \D y \, , \ \ \mbox{where} \ \ U_+ (y) =
\frac{\sinh (\mu |y|)} {\mu |y|} \, . \label{1+}
\end{equation}
On the other hand, the mean value property \eqref{MM+} is valid for $U_+$ over
$B_r$:
\begin{equation}
|B_r| \, a_3^+ (\mu r) = \int_{B_r} U_+ (y) \, \D y \, . \label{2+}
\end{equation}
 If we assume that $D \neq B_r$, then $G_i = D \setminus \overline{B_r}$ and $G_e
= B_r \setminus \overline D$ are bounded open sets such that $|G_e| = |G_i| \neq 0$,
which follows from the assumptions made about $D$ and $r$. Then, subtracting
\eqref{2+} from \eqref{1+}, we obtain
\begin{equation*}
0 = \int_{G_i} U_+ (y) \, \D y - \int_{G_e} U_+ (y) \, \D y > 0 \, .
\end{equation*}
Indeed, the difference is positive since $U_+ (y)$ (positive and monotonically
increasing with $|y|$) is greater than $[U_+ (y)]_{|y| = r}$ in $G_i$ and less than
$[U_+ (y)]_{|y| = r}$ in $G_e$, whereas $|G_i| = |G_e|$. The obtained contradiction
proves the result.
\end{proof}

\begin{remark}
{\rm The final part of this proof repeats literally the argument used in the proof
of Theorem~1; see \cite{Ku}, p.~947.}
\end{remark}

\begin{proof}[Proof of Theorem 4.]
In the acoustical case, we suppose, without loss of generality, that $D$ is located
so that $x_0$ coincides with the origin, and consider the following linear
combination of $E_\lambda^- (y, 0)$ and $E_\lambda^+ (y, 0)$:
\[ U_- (y) = \frac{\sin (\lambda |y|)} {\lambda |y|} \, ,
\]
As in the proof of Theorem~2, we arrive at the following consequence of
\eqref{MMEDH}
\begin{equation}
|D| \, a_3^- (\lambda r) = \int_D U_- (y) \, \D y \, , \label{1}
\end{equation}
where the condition $U_- (0) = 1$ is taken into account; cf. \eqref{1+}. Again,
assuming that $D \neq B_r$, we consider the bounded open sets $G_i = D \setminus
\overline{B_r}$ and $G_e = B_r \setminus \overline D$ such that $|G_e| = |G_i| \neq
0$; a consequence of the assumptions made about $D$ and $r$.

To obtain a contradiction we write the mean value property \eqref{MM-} for $U_-$
over $B_r$:
\begin{equation}
|B_r| \, a_3^- (\lambda r) = \int_{B_r} U_- (y) \, \D y \, . \label{2}
\end{equation}
Subtracting \eqref{2} from \eqref{1} and using the definition of $r$, we obtain
\begin{equation*}
0 = \int_{G_i} U_- (y) \, \D y - \int_{G_e} U_- (y) \, \D y < 0 \, .
\end{equation*}
Indeed, $U_- (y)$ monotonically decreases with $|y|$ in the whole $D$ because $D
\subset B_{r_0}$. Therefore, the difference is negative because $U_- (y)$ is
strictly greater than $[U_- (y)]_{|y| = r}$ in $G_e$ and strictly less than this
value in $G_i$, whereas $|G_i| = |G_e|$. The obtained contradiction proves the
theorem.
\end{proof} 

\begin{remark}
{\rm Here, the argument is similar to that in the proof of Theorem~2; indeed, both
rely on monotonicity of a certain solution to the corresponding equation. However,
there is an essential distinction between the two theorems concerning the size of a
domain. Indeed, no restriction on the size is imposed in Theorem~2. However, the
radially symmetric function $U_-$ decreases monotonically near the origin, but only
when $\lambda |y|$ belongs to a bounded interval adjacent to zero, for which reason
condition \eqref{exam} is imposed in Theorem~4.}
\end{remark}


{\small
\begin{thebibliography}{99}

\bibitem{ASZ} D. Aharonov, M.~M. Schiffer, L. Zalcman, ``Potato kugel'', {\sl Israel
J. Math.} {\bf 40} (1981), 331--339.

\bibitem{CL} G. Cupini, E. Lanconelli, ``On the harmonic characterization of domains
via mean value formulas'', {\sl Le Matematiche} {\bf 75} (2020), 331--352.

\bibitem{D} R.~J. Duffin, ``Yukawan potential theory'', {\sl J. Math. Anal. Appl.}
{\bf 35}, 105--130 (1971).

\bibitem{E} B. Epstein, ``On the mean-value property of harmonic functions'', {\sl
Proc. Amer. Math. Soc.} {\bf 13} (1962), 830.

\bibitem{ES} B. Epstein, M.~M. Schiffer, ``On the mean-value property of harmonic
functions'', {\sl J. Analyse Math.} {\bf 14} (1965), 109--111.

\bibitem{HN} W. Hansen, I. Netuka, ``Inverse mean value property of harmonic
functions'', {\sl Math. Ann.} {\bf 297} (1993), 147--156. Corrigendum: \textit{Math.
Ann.} {\bf 303} (1995), 373--375.

\bibitem{K} \"U, Kuran, ``On the mean value property of harmonic functions'', {\sl
Bull. London Math. Soc.} {\bf 4} (1972), 311--312.

\bibitem{Ku} N. Kuznetsov, ``Characterization of balls via solutions of the modified
Helmholtz equation,'' {\sl Comptes Rendus Math.} {\bf 359} (2021), 945--948.

\bibitem{Ku2} N. Kuznetsov, ``Mean value properties of solutions to the Helmholtz
and modified Helmholtz equations'', {\sl J. Math. Sci.} {\bf 257} (2021), 673--683.

\bibitem{Ku4} N. Kuznetsov, ``Inverse mean value properties (a survey)'', {\sl J.
Math. Sci.} {\bf 262} (2022), 275--290; see also preprint arxiv:2203.10601.

\bibitem{Ku1} N. Kuznetsov, ``Inverse mean value property of metaharmonic
functions'', {\sl J. Math. Sci.} (accepted); see also preprint arXiv:2203.14833.

\bibitem{L} E. Lanconelli, ````Potato kugel'' for sub-Laplacians'', {\sl Israel J.
Math.} {\bf 194} (2013), 277--283.

\bibitem{NV} I. Netuka, J. Vesel\'y, ``Mean value property and harmonic functions'',
{\sl Classical and Modern Potential Theory and Applications}. Kluwer, Dordrecht,
1994, pp. 359--398.

\bibitem{NC} C. Neumann, {\sl Allgemeine Untersuchungen \"uber das Newtonsche
Prinzip der Fernwirkungen}, Teubner, Leipzig, 1896.

\bibitem{S} M.~M. Schiffer, {\sl Selected Papers. Vol. 2}. P.~Duren, L.~Zalcman
(eds.) Springer, New York et al., 2014.

\bibitem{Wa} G.~N. Watson, {\sl A Treatise on the Theory of Bessel Functions}, 2nd
ed., Cambridge University Press, Cambridge, 1944.

\bibitem{Ve2} I. N. Vekua, {\sl New Methods for Solving Elliptic Equations}, North
Holland, Amsterdam, 1967.

\bibitem{Y} H. Yukawa, ``On the interaction of elementary particles,'' {\sl Proc.
Phys.-Math. Soc. Japan} {\bf 17} (1935), 48--57.

\end{thebibliography}
}
\end{document}